\documentclass[a4paper,12pt,oneside]{article} 
\usepackage{amssymb} 
\usepackage{amsmath} 
\usepackage{amsthm} 
\title{\textbf{Fractional dynamical systems and applications  in mechanics and economics}} 
\author{Mihai Boleantu \\  Dept. of Economics, West University 
of Timisoara, \\ email: mihaiboleantu01@yahoo.com 
\and Dumitru Opris \\ Dept. of Mathematics, 
West University of Timisoara, \\ email: miticaopris@yahoo.com} 

\begin{document}
\date{} %
\maketitle 
\begin{abstract} 
 
Using the fractional integration and differentiation on 
$\mathbb{R}$ we build the fractional jet fibre bundle 
on a differentiable manifold and we emphasize some 
important geometrical objects. Euler-Lagrange fractional 
equations are described. Some significant 
examples from mechanics and economics are presented. 
\end{abstract}  
\section{Introduction}  
\indent \par The operators of fractional differentiation 
have been introduced by Leibnitz, Liouville, Riemann, 
Grunwal and Letnikov \cite{Cresson}. The fractional 
derivatives and integrals are used in the description 
of some models in mechanics, physics \cite{Cresson}, 
economics \cite{Caputo} and 
medicine \cite{broasca}. The fractional variational 
calculus \cite{Agrawal}  
is an important instrument in the analysis of such 
models. The Euler-Lagrange equations are 
non-autonomous fractional differential equations 
in those models. 
\indent \par In this paper we present the fractional 
jet fibre bundle of order $k$ on a differentiable 
manifold as being $J^{\alpha k}({\mathbb R}, M)={\mathbb R}\times Osc^{\alpha k}(M)$, $\alpha \in (0,1)$, $k\in {\mathbb N}^{*}$. The fibre bundle $J^{\alpha k}$ is 
built in a similar way as the fibre bundle $E^{k}$ by 
R. Miron \cite{Miron}. Among the geometrical 
structures defined on $J^{\alpha}({\mathbb R}, M)$ 
we consider the dynamical fractional connection and 
the fractional Euler-Lagrange equations associated with 
a function defined on $J^{\alpha k}({\mathbb R}, M)$. \\     
\indent In section 2 we describe the fractional 
operators on ${\mathbb R}$ and some of their 
properties which are used in the paper. In section 3 
we describe the fractional osculator bundle 
of order $k$. In section 4 the fractional 
jet fibre bundle $J^{\alpha}({\mathbb R}, M)$ is 
defined, the fractional dynamical connection is 
built and the fractional Euler-Lagrange equations 
are established using the notion of fractional 
extremal value and classical extremal value on 
$J^{\alpha k}({\mathbb R}, M)$. In 
section 5 we consider some examples 
and applications. 
\section{Elements of fractional integration and 
differentiation on $\mathbb{R}$} 
\indent \par Let $f:[a, b]\to {\mathbb{R}}$ be 
an integrable function and $\alpha \in (0, 1)$. 
The left-sided (right-sided) fractional derivative 
of $f$ is the function 
\renewcommand{\arraystretch}{2} 
\begin{equation} \label{1)} \begin{array}{l} {({}_{-} D_{t}^{\alpha } f)(t)=\frac{1}{\Gamma (1-\alpha )} \frac{d}{dt} \int _{a}^{t}\frac{f(s)-f(a)}{(t-s)^{\alpha } } ds } \\ {({}_{+} D_{t}^{\alpha } f)(t)=\frac{1}{\Gamma (1-\alpha )} \frac{d}{dt} \int _{t}^{b}\frac{f(b)-f(s)}{(s-t)^{\alpha } } ds, } \end{array} \end{equation} 
where $t\in [a,\; b)$ and $\Gamma$ is Euler's 
gamma function. 
\newtheorem{theorem}{Theorem}  
\newtheorem{proposition}[theorem]{Proposition}  
\renewcommand{\arraystretch}{1.5} 
\begin{proposition} \normalfont(see \cite{Cresson})  
The operators ${}_{-} D_{t}^{\alpha }$ and 
${}_{+} D_{t}^{\alpha }$ have the properties:\\ 
1. If $f_{1}$ and $f_{2}$ are defined on 
$[a, b]$ and ${}_{-} D_{t}^{\alpha}$, 
${}_{+} D_{t}^{\alpha }$ exists, then 
\begin{equation} \label{2)} {}_{-} D_{t}^{\alpha } (c_{1} f_{1} +c_{2} f_{2} )(t)=c_{1} ({}_{-} D_{t}^
{\alpha } f_{1} )(t)+c_{2} ({}_{-} D_{t}^{\alpha } f_{2} )(t). \end{equation}  
2. If $\{ \alpha _{n} \} _{n\ge 0}$ is a real number 
sequence with 
$\mathop{\lim }\limits_{n\to \infty } \alpha _{n} =1$ 
then 
\begin{equation} \label{3)} \mathop{\lim }\limits_{n\to \infty } ({}_{-} D_{t}^{\alpha _{n} } f)(t)=({}_{-} D_{t}^{1} f)(t)=\frac{d}{dt} f(t). \end{equation} 
3. a) If $f(t)=c$, $t\in [a, b]$, $c\in {\mathbb{R}}$ 
then 
\begin{equation} \label{4)} ({}_{-} D_{t}^{\alpha } f)(t)=0.  \end{equation} 
b) If $f(t)=t^{\gamma }$, $t\in (a, b]$, 
$\gamma \in {\mathbb{R}}$, then 
\begin{equation} \label{5)} ({}_{-} D_{t}^{\alpha } f)(t)=\frac{t^{\gamma -\alpha } \Gamma (1+\gamma )}{\Gamma (1+\gamma -\alpha )} . \end{equation} 
c) If $f(t)=\frac{t^{\alpha } }{\Gamma (1+\alpha )}$, 
then 
\begin{equation} \label{6)} ({}_{-} D_{t}^{\alpha } f)(t)=1. \end{equation} 
4. If $f_{1}$ and $f_{2}$ are analytic functions on 
$[a, b]$ then 
\begin{equation} \label{7)} ({}_{-} D_{t}^{\alpha } (f_{1} f_{2} ))(t)=\sum _{k=0}^{\infty }\left(\begin{array}{c} {\alpha } \\ {k} \end{array}\right) ({}_{-} D_{t}^{\alpha -k} f_{1} )(t)\frac{d^{k} }{(dt)^{k} } f_{2} (t), \end{equation}  
where 
$\frac{d^{k} }{(dt)^{k} } =\frac{d}{dt} \circ \frac{d}{dt} \circ ...\circ \frac{d}{dt}$. \\ 
5. It also holds true 
\begin{equation} \label{8)} \int _{a}^{b}f_{1} (t)({}_{-} D_{t}^{\alpha } f_{2} )(t)dt=-\int _{a}^{b}f_{2} (t)({}_{+} D_{t}^{\alpha } f_{1} )(t)dt.   \end{equation} 
6. a) If $f:[a, b]\to {\mathbb R}$ admits fractional 
derivatives of order $a\alpha $, $a\in {\mathbb N}$, 
then 
\begin{equation} \label{9)} f(t+h)=E_{\alpha}((ht)^{\alpha } {}_{-} D_{t}^{\alpha } )f(t), \end{equation} 
where $E_{\alpha }$ is the Mittag-Leffler function 
given by 
\begin{equation} \label{10)} E_{\alpha } (t)=\sum _{a=0}^{\infty }\frac{t^{\alpha a} }{\Gamma (1+\alpha a)}. \end{equation} 
b) If $f:[a, b]\to {\mathbb R}$ is analytic and 
$0\in (a, b)$ then the fractional McLaurin series 
is 
\begin{equation} \label{ZEqnNum520508} f(t)=\sum _{a=0}^{\infty }\frac{t^{\alpha a} }{\Gamma (1+\alpha a)}  ({}_{-} D_{t}^{\alpha a} f)(t)\left|_{t=0} \right..  \end{equation} 
\end{proposition} 
\indent The physical and geometrical interpretation of 
the fractional derivative on ${\mathbb R}$ is suggested 
by the interpretation of the Stieltjes integral, because 
the integral used in the definition of the fractional 
derivative is a Riemann-Stieltjes integral \cite{Podlubny}.\\ 
\indent By definition, the left-sided (right-sided) 
fractional derivative of $f$, of order $\alpha$, 
$m=[\alpha]+1$, is the function 
\begin{equation} 
\label{eq12} 
\begin{array}{l} {D_{t}^{\alpha } f(t)=\frac{1}{\Gamma (m-\alpha )} \left(\frac{d}{dt} \right)^{m} \int _{-\infty }^{t}\frac{f(s)-f(0)}{(t-s)^{\alpha } } ds , \; 0\in(-\infty,t) } \\ {{}^{*} D_{t}^{\alpha } f(t)=\frac{1}{\Gamma (m-\alpha )} \left(-\frac{d}{dt} \right)^{m} \int _{t}^{\infty }\frac{f(s)-f(0)}{(s-t)^{\alpha } } ds , \; 0\in(t,\infty).} \end{array} \end{equation} 
If $\overline{{\rm supp}{\it f}}\subset [a,b]$, 
then $D_{t}^{\alpha}f={}_{-} D_{t}^{\alpha } f$, 
${}^{*} D_{t}^{\alpha } f={}_{+} D_{t}^{\alpha }f$.\\ 
\indent Let us consider the seminorms 
\begin{equation*} 
\begin{array}{l} {\left|x\right|_{J_{L}^{\alpha } ({\mathbb R})} =\left\| D_{t}^{\alpha } x\right\| _{L^{2} ({\mathbb R})} \; } \\ {\left|x\right|_{J_{{ R}}^{\alpha } ({\mathbb R})} =\left\| {}^{*} D_{t}^{\alpha } x\right\| _{L^{2} ({\mathbb R})} ,} \end{array} 
\end{equation*} 
and the norms 
\begin{equation*}  
\begin{array}{l} {\left\| x\right\| _{J_{L}^{\alpha } ({\mathbb R})} =\left(\left\| x\right\| _{L^{2} ({\mathbb R})}^{2} +\left|x\right|_{J_{L}^{\alpha } ({\mathbb R})}^{2} \right)^{1/2} \; } \\ {\left\| x\right\| _{J_{{R}}^{\alpha } ({\mathbb R})} =\left(\left\| x\right\| _{L^{2} ({\mathbb R})}^{2} +\left|x\right|_{J_{{R}}^{\alpha } ({\mathbb R})}^{2} \right)^{1/2},} \end{array} 
\end{equation*} 
and $J_{_{0}L}^{\alpha } ({\mathbb R})$, 
$J_{_{0} R}^{\alpha } ({\mathbb R})$ 
the closures of $C_{0}^{\infty} ({\mathbb R})$ 
with respect to the two norms from above, respectively. 
In \cite{Cresson} it is proved that the 
operators $D_{t}^{\alpha }$ and ${}^{*} D_{t}^{\alpha}$ 
satisfy the properties: 
\begin{proposition} 
Let $I\subset {\mathbb R}$ and let 
$J_{_{0}L }^{\alpha } (I)$ 
and $J_{_{0} R} ^{\alpha } (I)$ be the closures 
of $C_{0}^{\infty } (I)$ with respect to the norms 
from above. For any 
$x\in J_{_{0}L }^{\beta } (I)$, $0<\alpha <\beta $, 
the following relation holds: 
\begin{equation*} 
D_{t}^{\beta } x(t)=D_{t}^{\alpha } D_{t}^{\beta -\alpha } x(t). \end{equation*} 
For any $x\in J_{_{0} R }^{\beta } (I)$, 
$0<\alpha <\beta $, it also holds 
\begin{equation*} 
{}^{*} D_{t}^{\beta } x(t)={}^{*} D_{t}^{\alpha } {}^{*} D_{t}^{\beta -\alpha } x(t). \end{equation*}  
\end{proposition}  
In the following we shall consider the fractional 
derivatives defined above. 
\section{The fractional osculator bundle of order $k$ on a differentiable manifold} 
\indent \par Let $\alpha \in (0, 1]$ be fixed and 
$M$ a differentiable manifold of dimension $n$. 
Two curves $\rho ,\; \sigma :I\to {\mathbb R}$, 
with $\rho(0) = \sigma(0)=x_{0} \in M$, 
$0\in I$, have a fractional contact $\alpha$ 
of order $k\in {\mathbb N}^{*}$ in $x_{0} $, 
if for any $f\in {\mathcal F}(U)$, $x_{0} \in U$, 
$U$ a chart on $M$, it holds 
\begin{equation} \label{ZEqnNum286695} D_{t}^{\alpha a} (f\circ \rho )\left|_{t=0} \right. =D_{t}^{\alpha a} (f\circ \sigma )\left|_{t=0} \right.  \end{equation} 
where $a=\overline{1,\; k}$. The relation 
\eqref{ZEqnNum286695} is an equivalence relation. 
The equivalence class $[\rho ]_{x_{0} }^{\alpha k}$ 
is called the fractional \textit{k}-osculator space 
of $M$ in $x_{0}$ and it will be denoted by 
$Osc_{x_{0} }^{\alpha k} (M)$. 
If the curve $\rho :I\to M$ is given by 
$x^{i} =x^{i} (t)$, $t\in I$, $i=\overline{1,\; n}$, 
then, considering the formula \eqref{ZEqnNum520508}, 
the class $[\rho ]_{x_{0} }^{\alpha k}$, 
may be written as 
\begin{equation} \label{14)} x^{i} (t)=x^{i} (0)+\frac{t^{\alpha } }{\Gamma (1+\alpha )} D_{t}^{\alpha } x^{i} (t)\left|_{t=0} \right. +...+\frac{t^{\alpha k} }{\Gamma (1+\alpha k)} D_{t}^{\alpha k} x^{i} (t)\left|_{t=0} \right. , \end{equation} 
where $t\in (-\varepsilon , \varepsilon )$. 
We shall use the notation 
\begin{equation} \label{ZEqnNum283756} x^{i} (0)=x^{i},\quad  y^{i(\alpha a)} =\frac{1}{\Gamma (1+\alpha a)} D_{t}^{\alpha a} x^{i} (t)\left|_{t=0} \right. , \end{equation} 
for $i=\overline{1, n}$ and $a=\overline{1, k}$.\\ 
By definition, the fractional osculator bundle of order 
$r$ is the fibre bundle  
$(Osc^{\alpha k}(M), M)$ where 
$Osc^{\alpha k} (M)=\bigcup _{x_{0} \in M}Osc_{x_{0} }^{\alpha k}(M)$ 
and \\ 
$\pi _{0}^{\alpha k}:Osc^{\alpha k} (M)\to M$ 
is defined by 
$\pi_{0} ^{\alpha k} ([\rho ]_{x_{0} }^{\alpha k} )=x_{0}$, 
$(\forall )[\rho ]_{x_{0} }^{\alpha k} \in Osc^{\alpha k} (M)$. \\  
\indent For $f\in {\mathcal F}(U)$, the fractional derivative 
of order $\alpha $, $\alpha \in (0,\; 1)$, with respect to 
the variable $x^{i} $, is defined by 
\begin{equation} \label{16)} \begin{array}{l} {(D_{x^{i} }^{\alpha } f)(x)=} \\ {\frac{1}{\Gamma (1-\alpha )} \frac{\partial }{\partial x^{i} } \int _{a^{i} }^{x^{i} }\frac{f(x^{1} ,...,x^{i-1} ,s,x^{i+1} ,...,x^{n} )-f(x^{1} ,...,x^{i-1} ,a^{i} ,x^{i+1} ,...,x^{n} )}{(x^{i} -s)^{\alpha } }  ds,} \end{array} \end{equation} 
where $x^{i}$ are the coordinate functions on $U$, 
$\frac{\partial }{\partial x^{i} } $, 
$i=\overline{1, n}$, is the canonical base of the 
vector fields on $U$ and 
$U_{ab} =\{ x\in U,\; a^{i} \le x^{i} \le b^{i} ,\; i=\overline{1,n}\} \subset U$. \\ 
Let $U,\; U'\subset M$ be two charts on $M$, 
$U\cap U'\ne \emptyset$ and consider the change of variable 
\begin{equation} \label{ZEqnNum197194} \bar{x}^{i} =\bar{x}^{i} (x^{1} ,...,x^{n} ) \end{equation} 
with $\det \left(\frac{\partial \bar{x}^{i} }{\partial x^{j} } \right)\ne 0$. 
Let 
$\left\{dx^{i} \right\}_{i=\overline{1,\; n}}$ 
be the canonical base of 1-forms of 
${\mathcal D}^{1} (U)$ and let us define the 
1-forms 
$d(x^{i} )^{\alpha } =\alpha (x^{i} )^{\alpha -1} dx^{i} $, 
$i=\overline{1, n}$. 
The exterior differential 
$d^{\alpha } :{\mathcal F}(U\cap U')\to {\mathcal D}^{1} (U\cap U')$ 
is defined by 
\begin{equation} \label{ZEqnNum198901} d^{\alpha } =d(x^{j} )^{\alpha }D_{x^{j} }^{\alpha }=d(\bar{x}^{j} )^{\alpha} D_{\bar{x}^{j} }^{\alpha } . \end{equation} 
Using \eqref{ZEqnNum198901} and the property 
$D_{x^{i} }^{\alpha } \left(\frac{(x^{i} )^{\alpha } }{\Gamma (1+\alpha )} \right)=1$, 
it follows that 
\begin{equation} \label{ZEqnNum965179} d(x^{j} )^{\alpha } =\frac{1}{\Gamma (1+\alpha )} D_{\bar{x}^{i} }^{\alpha } (x^{j} )^{\alpha } d(\bar{x}^{i} )^{\alpha }.  \end{equation} 
Using the notation 
\begin{equation} \label{20)} \mathop{J_{i}^{j} }\limits^{\alpha } (x,\; \bar{x})=\frac{1}{\Gamma (1+\alpha )} D_{\bar{x}^{i} }^{\alpha } (x^{j} )^{\alpha }, \end{equation} 
from \eqref{ZEqnNum965179} we get 
\begin{equation} \label{ZEqnNum547724} d(x^{j} )^{\alpha } =\mathop{J_{i}^{j} }\limits^{\alpha } (x,\; \bar{x})d(\bar{x}^{i} )^{\alpha }.  \end{equation} 
From \eqref{ZEqnNum547724} it follows that 
\begin{equation} \label{22)} \mathop{J_{i}^{j} }\limits^{\alpha } (x,\; \bar{x})\mathop{J_{h}^{i} }\limits^{\alpha } (x,\; \bar{x})=\delta _{h}^{j} . \end{equation} 
Consider $x^{i} =x^{i} (t)$ and 
$\bar{x}^{i} (t)=\bar{x}^{i} (x(t))$, 
$i=\overline{1, n}$, $t\in I$. 
Applying the operator  $D_{t}^{\alpha }$ 
we get 
\begin{equation} \label{23)} (D_{t}^{\alpha } \bar{x}^{i} )(t)=D_{x^{j} }^{\alpha } \bar{x}^{i} (x)(D_{t}^{\alpha } x^{j} )(t)=\mathop{J_{j}^{i} }\limits^{\alpha } (\bar{x},\; x)(D_{t}^{\alpha } x^{j} )(t). \end{equation} 
Considering the notation from \eqref{ZEqnNum283756} 
we have 
\begin{equation} \label{ZEqnNum957067} y^{i(\alpha )} =\mathop{J_{j}^{i} }\limits^{\alpha } (\bar{x},\; x)\bar{y}^{j(\alpha )} . \end{equation} 
Also, from \eqref{ZEqnNum283756} we deduce 
\begin{equation} \label{25)} D_{t}^{\alpha } y^{i(\alpha a)} =\frac{\Gamma (\alpha a)}{\Gamma (\alpha (a-1))} y^{i(\alpha a)} , \end{equation} 
where $i=\overline{1, n}$. Applying the operator 
$D_{t}^{\alpha }$ in the relation 
\eqref{ZEqnNum957067} we find 
\begin{equation} \label{ZEqnNum440443} \begin{array}{l} {\frac{\Gamma (\alpha (a-1))}{\Gamma (\alpha )} \bar{y}^{i(\alpha a)} =\Gamma (1+\alpha )\mathop{J_{j}^{i} }\limits^{\alpha } (\bar{y}^{\alpha (a-1)} ,\; x)y^{j(\alpha )} +} \\ {\frac{\Gamma (2\alpha )}{\Gamma (\alpha )} \mathop{J_{j}^{i} }\limits^{\alpha } (y^{(\alpha (a-1))} ,\; y^{\alpha } )y^{j(2\alpha )} +...+\frac{\Gamma (2\alpha )}{\Gamma (\alpha )} \mathop{J_{j}^{i} }\limits^{\alpha } (\bar{y}^{\alpha (a-1)} ,\; y^{\alpha b} )y^{j((b+1)\alpha )} +} \\ {...+\frac{\Gamma (\alpha (a-1))}{\Gamma (\alpha )} y^{i(\alpha a)} ,} \end{array} \end{equation} 
where $a=\overline{1, k}$. 
\begin{proposition} \normalfont(see \cite{Albu}, \cite{Cottrill}) \\   
a) The coordinate transformation on $Osc^{(\alpha k)}(M)$, 
\\ 
$(x^{i} ,y^{i(\alpha )}  ,...,y^{i(\alpha k)} )\to (\bar{x}^{i} ,\bar{y}^{i(\alpha )}  ,...,\bar{y}^{i(\alpha k)} )$ 
are given by the formulas \eqref{ZEqnNum197194} 
and \eqref{ZEqnNum440443}. \\ 
b) The operators $D_{x^{i} }^{\alpha }$ and the 
1-forms  $(dx^{i} )^{\alpha } $, $i=\overline{1, n}$, 
transform by the formulas 
\begin{equation} \label{27)} \begin{array}{l} {D_{\bar{x}^{i} }^{\alpha } =\mathop{J_{j}^{i} }\limits^{\alpha } (x,\; \bar{x})D_{x^{j} }^{\alpha } } \\ {d(\bar{x}^{i})^{\alpha} =\mathop{J_{j}^{i} }\limits^{\alpha } (\bar{x},\; x)d(x^{j} )^{\alpha } .} \end{array} \end{equation} 
\end{proposition} 
\section{The fractional jet bundle of order $k$\\ 
on a differentiable manifold; geometrical 
objects } 
\indent \par By definition, the $k$-order fractional 
jet bundle is the space 
$J^{\alpha k} ({\mathbb R},\; M)={\mathbb R}\times Osc^{k\alpha } (M).$ 
A system of local coordinates on 
$J^{\alpha k} ({\mathbb R},\; M)$ will be denoted by 
$(t,x,y^{(\alpha )} ,y^{(2\alpha )} ,...,y^{(k\alpha )} )$. 
Consider the projections 
$\pi _{0}^{\alpha k} :J^{\alpha k}({\mathbb R},\; M)\to M$ 
defined by 
\begin{equation} \label{28)} \begin{array}{l} {\pi _{0}^{\alpha k} (t,x,y^{(\alpha )} ,...,y^{(\alpha k)} )=x.} \end{array} \end{equation} 
Let $U,\; U'\subset M$ be two charts on $M$ with 
$U\cap U'\ne \emptyset $, 
$(\pi _{0}^{\alpha } )^{-1} (U),\; (\pi _{0}^{\alpha } )^{-1} (U')\subset J^{\alpha } ({\mathbb R},\; M)$ 
the corresponding charts on 
$J^{\alpha } ({\mathbb R},\; M)$ 
and, respectively, the corresponding 
coordinates 
$(x^{i} )$, $(\bar{x}^{i} )$ and 
$(t,x^{i} ,y^{i(\alpha )} )$,  
$(t,\bar{x}^{i} ,\bar{y}^{i(\alpha )} )$. 
The transformations of coordinates are given by 
\begin{equation} \label{ZEqnNum534725} \begin{array}{l} {\bar{x}^{i} =\bar{x}^{i} (x^{1} ,...,x^{n} )} \\ {\bar{y}^{i(\alpha )} =\mathop{J}\limits^{\alpha } (x,\bar{x})y^{i(\alpha )} .} \end{array} \end{equation} 
\indent Consider the functions  
$(t)^{\alpha } $, $(x^{i} )^{\alpha }$, 
$(y^{i(\alpha )} )^{\alpha } \in {\mathcal F}((\pi_{0} ^{\alpha } )^{-1} (U))$, 
the 1-forms 
$\frac{1}{\Gamma (1+\alpha )} d(t)^{\alpha } $, 
$\frac{1}{\Gamma (1+\alpha )} d(x^{i} )^{\alpha } $, 
$\frac{1}{\Gamma (1+\alpha )} d(y^{i(\alpha )} )^{\alpha } \in {\mathcal D}^{1} ((\pi_{0} ^{\alpha } )^{-1} (U))$ 
and the operators 
$D_{t}^{\alpha } $, $D_{x^{i} }^{\alpha } $, 
$D_{y^{i(\alpha )} }^{\alpha } $ on 
$(\pi_{0} ^{\alpha } )^{-1} (U)$, $i=\overline{1, n}$. 
The following relations hold: 
\begin{equation} \label{ZEqnNum227932} \begin{array}{l} {D_{t}^{\alpha } (\frac{1}{\Gamma (1+\alpha )} t^{\alpha } )=1,\; \; D_{x^{i} }^{\alpha } (\frac{1}{\Gamma (1+\alpha )} (x^{j} )^{\alpha } )=\delta _{i}^{j} ,} \\ {D_{y^{i(\alpha )} }^{\alpha } (\frac{1}{\Gamma (1+\alpha )} (y^{j(\alpha )} )^{\alpha } )=\delta _{i}^{j} ,\; \; \frac{1}{\Gamma (1+\alpha )} d(t^{\alpha } )(D_{t}^{\alpha } )=1,} \\ {\frac{1}{\Gamma (1+\alpha )} d(x^{i} )^{\alpha } (D_{x^{j} }^{\alpha } )=\delta _{j}^{i} ,\; \; \frac{1}{\Gamma (1+\alpha )} d(y^{i(\alpha )} )^{\alpha } (D_{y^{j(\alpha )} }^{\alpha } )=\delta _{j}^{i} .} \end{array} \end{equation} 
On $J^{\alpha } ({\mathbb R},\; M)$ we may define 
the canonical structures 
\begin{equation} \label{ZEqnNum331309} \begin{array}{l} {\mathop{\theta_{1} }\limits^{\alpha } =d(t^{\alpha } )\otimes (D_{t}^{\alpha } +y^{i(\alpha )} D_{x^{i} }^{\alpha } )} \\ {\mathop{\theta_{2} }\limits^{\alpha }  =\mathop{\theta ^{i} }\limits^{\alpha } \otimes D_{x^{i} }^{\alpha } ,\; \; \mathop{\theta ^{i} }\limits^{\alpha } =\frac{1}{\Gamma (1+\alpha )} (d(x^{i} )^{\alpha } -y^{i(\alpha )} d(t)^{\alpha } )} \\ {\mathop{S}\limits^{\alpha } =\mathop{\theta ^{i} }\limits^{\alpha } \otimes D_{y^{i(\alpha )} }^{\alpha } } \\ {\mathop{V_{i} }\limits^{\alpha } =D_{y^{i(\alpha )} }^{\alpha } .} \end{array} \end{equation} 
Using \eqref{ZEqnNum534725} it is easy to show that 
the structures \eqref{ZEqnNum331309} 
have geometrical character. The space of the operators 
generated by the operators 
$\{ D_{t}^{\alpha } ,D_{x^{i} }^{\alpha } ,D_{y^{i(\alpha )} }^{\alpha } \} $, $i=\overline{1,\; n}$, 
will be denoted by 
$\chi ^{\alpha } ((\pi _{0}^{\alpha } )^{-1} (U))$. 
For $\alpha \to 1$ the space of these operators 
represents the space of the vector fields 
on $\pi _{0}^{-1} (U)$. \\ 
\indent A vector field 
$\mathop{\Gamma }\limits^{\alpha } \in \chi ^{\alpha } ((\pi _{0}^{\alpha } )^{-1} (U))$ 
is called \textit{FODE} (fractional ordinary differential 
equation) iff 
\begin{equation} \label{32)} \begin{array}{l} {d(t)^{\alpha } (\mathop{\Gamma }\limits^{\alpha } )=1} \\ {\mathop{\theta ^{i} }\limits^{\alpha } (\mathop{\Gamma }\limits^{\alpha } )=0,} \end{array} \end{equation}  
for $i=\overline{1, n}$. In local coordinates 
\textit{FODE} is given by 
\begin{equation} \label{33)} \mathop{\Gamma }\limits^{\alpha } =D_{t}^{\alpha } +y^{i(\alpha )} D_{x^{i} }^{\alpha } +F^{i} D_{y^{i(\alpha )} }^{\alpha } , \end{equation} 
where 
$F^{i} \in C^{\infty } ((\pi _{0}^{\alpha } )^{-1} (U))$, 
$i=\overline{1, n}$. 
The integral curves of the field \textit{FODE} 
are the solutions of the fractional differential 
equation (\textit{EDF}) 
\begin{equation} \label{34)} D_{t}^{2\alpha } x^{i} (t)=F^{i} (t,x(t),D_{t}^{\alpha } x(t)), \qquad i=\overline{1, n}. \end{equation}  
\indent The fractional dynamical connection on 
$J^{\alpha } ({\mathbb R},\; M)$ is defined by 
the fractional tensor fields 
$\mathop{H}\limits^{\alpha }$ of type $(1,1)$ which satisfy the conditions 
\begin{equation} \label{35)} \begin{array}{l} {\mathop{\theta _{1} }\limits^{\alpha } \circ \mathop{H}\limits^{\alpha } =0\; } \\ {\mathop {\theta _{2}}\limits^{\alpha} \circ \mathop{H}\limits^{\alpha }} ={\mathop{\theta _{2} }\limits^{\alpha } \; } \\ {\mathop{H}\limits^{\alpha }\left|_{\mathop{V}\limits^{\alpha } } \right. =-id\left|_{\mathop{V}\limits^{\alpha } } \right. ,\; } \end{array} \end{equation} 
where $\mathop{V}\limits^{\alpha }$ is formed by 
operators generated by 
$\{ D_{y^{i(\alpha )} }^{\alpha } \} _{i=\overline{1,n}}$. 
In the chart 
$(\pi _{0}^{\alpha } )^{-1} (U)$ the fractional tensor field 
$\mathop{H}\limits^{\alpha }$ has the expression 
\begin{equation} \label{36)} \begin{array}{l} {\mathop{H}\limits^{\alpha } =(\mathop{H}\limits^{1} d(t)^{\alpha } +\mathop{H_{j}^{i} }\limits^{2} d(x^{i} )^{\alpha } +\mathop{H_{i} }\limits^{3} d(y^{i(\alpha )} )^{\alpha } )\otimes D_{t}^{\alpha } +} \\ {(\mathop{H_{j}^{i} }\limits^{4} (dt)^{\alpha } +\mathop{H_{j}^{i} }\limits^{5} d(x^{j} )^{\alpha } +\mathop{H_{j}^{i} }\limits^{6} d(y^{i(\alpha )} )^{\alpha } )\otimes D_{x^{i} }^{\alpha } +} \\ {(\mathop{H^{i} }\limits^{7} d(t)^{\alpha } +\mathop{H_{j}^{i} }\limits^{8} d(x^{j} )^{\alpha } +\mathop{H_{j}^{i} }\limits^{9} d(y^{i(\alpha )} )^{\alpha } )\otimes D_{y^{i(\alpha )} }^{\alpha } .} \end{array} \end{equation} 
The tensor field $\mathop{H}\limits^{\alpha } $ has 
a geometrical character, fact which results by 
using the relations 
\eqref{ZEqnNum534725}, and is called a 
$d^{\alpha }$-tensor field. Using the relations 
\eqref{ZEqnNum227932} and \eqref{ZEqnNum331309} we get 
\begin{proposition} 
a) The fractional dynamical connection 
$\mathop{H}\limits^{\alpha }$, in the chart 
$(\pi _{0}^{\alpha } )^{-1} (U)$, 
is given by 
\begin{equation} \label{37)} \begin{array}{l} {\mathop{H}\limits^{\alpha } =\frac{1}{\Gamma (1+\alpha )} [(-y^{i(\alpha )} D_{x^{i} }^{\alpha } +H^{i} D_{y^{i(\alpha )} }^{\alpha } )\otimes d(t)^{\alpha } +} \\ {(D_{x^{i} }^{\alpha } +H_{i}^{j} D_{y^{j(\alpha )} }^{\alpha } )\otimes d(x^{i} )^{\alpha } -D_{y^{i(\alpha )} }^{\alpha }\otimes  d(y^{i(\alpha )} )^{\alpha } ].} \end{array} \end{equation} 
b) The fractional dynamical connection 
$\mathop{H}\limits^{\alpha } $ 
defines a $f(3,-1)$ fractional structure on 
$J^{\alpha } ({\mathbb R},\; M)$, i.e., 
$\left(\mathop{H}\limits^{\alpha } \right)^{3} =\mathop{H}\limits^{\alpha } $. \\ 
c) The fractional tensor fields 
$\mathop{l}\limits^{\alpha } $ 
and $\mathop{m}\limits^{\alpha } $ 
which are defined by 
\begin{equation} \label{38)} \begin{array}{l} {\mathop{l}\limits^{\alpha } =\mathop{H}\limits^{\alpha } \circ \mathop{H}\limits^{\alpha } } \\ {\mathop{m}\limits^{\alpha } =-\mathop{H}\limits^{\alpha } \circ \mathop{H}\limits^{\alpha } +I,} \end{array} \end{equation} 
where $I$ is the identity map, satisfy the relations 
\begin{equation} \label{39)} \begin{array}{l} {\mathop{l}\limits^{\alpha } \circ \mathop{l}\limits^{\alpha } =\mathop{l}\limits^{\alpha } ,\; \; \mathop{m}\limits^{\alpha } \circ \mathop{m}\limits^{\alpha } =\mathop{m}\limits^{\alpha } \circ \mathop{l}\limits^{\alpha } ,\; \; \mathop{l}\limits^{\alpha } +\mathop{m}\limits^{\alpha } =I} \\ {\mathop{l}\limits^{\alpha } (D_{t}^{\alpha } )=-y^{i(\alpha )} D_{x^{i} } -(y^{i(\alpha )} \mathop{H_{i}^{j} }\limits^{\alpha } +\mathop{H^{j} }\limits^{\alpha } )D_{y^{i(\alpha )} } } \\ {\mathop{l}\limits^{\alpha } (D_{x^{i} }^{\alpha } )=D_{x^{i} }^{\alpha } ,\; \; \mathop{l}\limits^{\alpha } (D_{y^{i(\alpha )} }^{\alpha } )=D_{y^{i(\alpha )} }^{\alpha } } \\ {\mathop{m}\limits^{\alpha } (D_{t}^{\alpha } )=D_{t}^{\alpha } +y^{i(\alpha )} D_{x^{i} }^{\alpha } +(y^{i(\alpha )} \mathop{H_{i}^{j} }\limits^{\alpha } +\mathop{H^{j} }\limits^{\alpha } )D_{y^{i(\alpha )} }^{\alpha } } \\ {\mathop{m}\limits^{\alpha } (D_{x^{i} }^{\alpha } )=0,\; \; \mathop{m}\limits^{\alpha } (D_{y^{i(\alpha )} }^{\alpha } )=0.} \end{array} \end{equation} 
d) The fractional vector field 
$\mathop{\Gamma }\limits^{\alpha } \in \chi ^{\alpha } (J^{\alpha } ({\mathbb R},\; M))$ 
given by 
\begin{equation} \label{40)} \mathop{\Gamma }\limits^{\alpha }  =\mathop{m}\limits^{\alpha}(D_{t}^{\alpha } )=D_{t}^{\alpha } +y^{i(\alpha )} D_{x^{i} }^{\alpha } +(y^{i(\alpha )} \mathop{H_{i}^{j} }\limits^{\alpha } +\mathop{H^{j} }\limits^{\alpha } )D_{y^{j(\alpha )} }^{\alpha }  \end{equation} 
defines a field \textit{FODE} associated to 
the fractional dynamical connection. 
The integral curves are the solutions of the 
\textit{EDF} 
\begin{equation} \label{41)} D_{t}^{2\alpha } x^{i} (t)=D_{t}^{\alpha } x^{i} (t)\mathop{H_{i}^{j} }\limits^{\alpha } +\Gamma (1+\alpha )\mathop{H^{j} }\limits^{\alpha }  \end{equation} 
where $\mathop{H_{i}^{j} }\limits^{\alpha } $ and 
$\mathop{H^{j} }\limits^{\alpha } $ 
are functions of 
$(t,\; x(t),\; y^{(\alpha) } (t))$. 
\end{proposition} 
Let 
$L\in C^{\infty } (J^{\alpha } ({\mathbb R},\; M))$ 
be a fractional Lagrange function. 
By definition, the Cartan fractional 1-form is 
the 1-form $\mathop{\theta _{L} }\limits^{\alpha }$ 
given by 
\begin{equation} \label{42)} \mathop{\theta _{L} }\limits^{\alpha } =Ld(t)^{\alpha } +\mathop{S}\limits^{\alpha } (L). \end{equation} 
We call the Cartan fractional 2-form, the 2-form 
$\mathop{\omega _{L} }\limits^{\alpha }$ 
given by 
\begin{equation} \label{43)} \mathop{\omega _{L} }\limits^{\alpha } =d^{\alpha } \mathop{\theta _{L} }\limits^{\alpha }  \end{equation} 
where $d^{\alpha } $ is the 
fractional exterior differential: 
\begin{equation} \label{44)} d^{\alpha } =d(t)^{\alpha } D_{t}^{\alpha } +d(x^{i} )^{\alpha } D_{x^{i} }^{\alpha } +d(y^{i(\alpha )} )^{\alpha } D_{y^{i(\alpha )} }^{\alpha } . \end{equation} 
In the chart $(\pi _{0}^{\alpha } )^{-1} (U)$, 
$\mathop{\theta _{L} }\limits^{\alpha } $ and 
$\mathop{\omega _{L} }\limits^{\alpha } $ 
are given by 
\begin{equation} \label{45)} \begin{array}{l} {\mathop{\theta _{L} }\limits^{\alpha } =(L-\frac{1}{\Gamma (1+\alpha )} y^{i(\alpha )} D_{y^{i(\alpha )} }^{\alpha } (L)d(t)^{\alpha } +\frac{1}{\Gamma (1+\alpha )} D_{y^{i(\alpha )} }^{\alpha } (L)d(x^{i} )^{\alpha } } \\ {\mathop{\omega _{L} }\limits^{\alpha } =A_{i} d(t)^{\alpha } \wedge d(x^{i} )^{\alpha } +B_{i} d(t^{\alpha } )\wedge d(y^{i(\alpha )} )^{\alpha } +} \\ {\; \; \; \; \; \; \; A_{ij} d(x^{i} )^{\alpha } \wedge d(x^{j} )^{\alpha } +B_{ij} d(x^{i} )^{\alpha } \wedge d(y^{j(\alpha )} )^{\alpha } ,} \end{array} \end{equation} 
where 
\begin{equation} \label{46)} \begin{array}{l} {A_{i} =\frac{1}{\Gamma (1+\alpha )} D_{t}^{\alpha } D_{y^{i(\alpha )} }^{\alpha } (L)+\frac{1}{\Gamma (1+\alpha )} y^{j(\alpha )} D_{x^{i} }^{\alpha } D_{y^{j(\alpha )} }^{\alpha } (L)-D_{x^{i} }^{\alpha } (L)} \\ {B_{i} =\frac{1}{\Gamma (1+\alpha )} D_{y^{i(\alpha )} }^{\alpha } (y^{j(\alpha )} D_{j}^{\alpha } (\alpha )(L))} \\ {A_{ij} =D_{x^{i} }^{\alpha } D_{y^{i(\alpha )} }^{\alpha } (L),\; \; B_{ij} =-D_{y^{j(\alpha )} }^{\alpha } D_{y^{i(\alpha )} }^{\alpha } (L).} \end{array} \end{equation} 
\begin{proposition} 
If $L$ is regular (i.e., 
$\det \left(\frac{\partial ^{2} L}{\partial y^{i(\alpha )} \partial y^{j(\alpha )} } \right)\ne 0$) 
then there exists a fractional field 
\textit{FODE} 
$\mathop{\Gamma _{L} }\limits^{\alpha }$ 
such that 
$i_{\mathop{\Gamma _{L} }\limits^{\alpha } } \mathop{\omega _{L} }\limits^{\alpha } =0$. 
In the chart $(\pi _{0}^{\alpha } )^{-1} (U)$ 
we have 
\begin{equation} \label{47)} \mathop{\Gamma _{L} }\limits^{\alpha } =D_{t}^{\alpha } +y^{i(\alpha )} D_{x^{i} } +\mathop{M^{i} }\limits^{\alpha } D_{y^{i(\alpha )} }^{\alpha } , \end{equation} 
where 
\begin{equation} \label{48)} \begin{array}{l} {\mathop{M^{i} }\limits^{\alpha } =g^{ik} (D_{k}^{\alpha } (L)-d_{t}^{\alpha } (\frac{\partial ^{\alpha } L}{\partial y^{k(\alpha )} } )} \\ {d_{t}^{\alpha } =D_{t}^{\alpha } +y^{i(\alpha )} D_{x^{i} }^{\alpha } } \\ {(g^{ik} )=(D_{y^{i(\alpha )} }^{\alpha } D_{y^{k(\alpha )} }^{\alpha } (L))^{-1} .} \end{array} \end{equation} 
\end{proposition} 
An important structure on 
$J^{\alpha } ({\mathbb R},\; M)$ 
is described by the fractional Euler-Lagrange equations. 
Let $c:t\in [0,1]\to (x^{i} (t))\in M$ be a 
parameterized curve, such that 
$Imc\subset U\subset M$. 
The extension of the curve $c$ to 
$J^{\alpha } ({\mathbb R},\; M)$ 
is the curve 
$c^{\alpha } :t\in [0,1]\to (t,x^{i} (t),y^{i(\alpha )} (t))\in J^{\alpha } ({\mathbb R},\; M)$. 
Consider 
$L\in C^{\infty } (J^{\alpha } ({\mathbb R},\; M))$. 
The action of $L$ along the curve $c^{\alpha }$ 
is defined by 
\begin{equation} \label{49)} {\mathcal A}(c^{\alpha } )=\int _{0}^{1}L(t,x(t),y^{\alpha } (t))dt . \end{equation} 
Let 
$c_{\varepsilon } :t\in [0,1]\to (x^{i} (t,\varepsilon ))\in M$ 
be a family of curves, where $\varepsilon $ is 
sufficiently small so that 
$Imc_{\varepsilon } \subset U$, $c_{0} (t)=c(t)$, 
$D_{\varepsilon }^{\alpha } c_{\varepsilon }(0)=D_{\varepsilon }^{\alpha } c_{\varepsilon}(1)=0$. 
The action of $L$ along the curves 
$c_{\varepsilon } $ is 
\begin{equation} \label{ZEqnNum353506} {\mathcal A}(c_{\varepsilon }^{\alpha } )=\int _{0}^{1}L(t,x(t,\varepsilon ),y^{\alpha } (t,\varepsilon ))dt , \end{equation} 
where 
$y^{i(\alpha )} (t,\varepsilon )=\frac{1}{\Gamma (1+\alpha )} D_{t}^{\alpha } x^{i} (t,\varepsilon )$. 
The action \eqref{ZEqnNum353506} 
has a fractional extremal value if 
\begin{equation} \label{51)} D_{\varepsilon }^{\alpha } {\mathcal A}(c_{\varepsilon }^{\alpha } )\left|_{\varepsilon =0} \right. =0. \end{equation} 
The action \eqref{ZEqnNum353506} 
has an extremal value if 
\begin{equation} \label{52)} D_{\varepsilon }^{1} {\mathcal A}(c_{\varepsilon }^{\alpha } )\left|_{\varepsilon =0} \right. =0. \end{equation} 
Using the properties of the fractional derivative we obtain 
\begin{proposition} 
a) A necessary condition for the action 
\eqref{ZEqnNum353506} to reach a fractional extremal value 
is that $c(t)$ satisfies the fractional Euler-Lagrange 
equations 
\begin{equation} \label{ZEqnNum157890} \begin{array}{l} {D_{x^{i} }^{\alpha } L-d_{t}^{2\alpha } (D_{y^{i(\alpha )} }^{\alpha } L)=0} \\ {d_{t}^{\alpha } =D_{t}^{\alpha } +y^{i(\alpha )} D_{x^{i} }^{\alpha } +y^{i(2\alpha )} D_{y^{i(\alpha )} }^{\alpha } ,} \end{array} \end{equation} 
where $i=\overline{1, n}$. \\ 
b) A necessary condition for the action 
\eqref{ZEqnNum353506} to reach an extremal value is that 
$c(t)$ satisfies the Euler-Lagrange equations 
\begin{equation} \label{ZEqnNum578225} \begin{array}{l} {D_{x^{i} }^{1} L-d_{t}^{2} (D_{y^{i(\alpha )} }^{1} L)=0} \\ {d_{t}^{2} =D_{t}^{\alpha } +y^{i(\alpha )} D_{x^{i} }^{1} +y^{i(2\alpha )} D_{y^{i(\alpha )} }^{1} ,} \end{array} \end{equation} 
where $i=\overline{1, n}$. 
\end{proposition} 
The equations \eqref{ZEqnNum157890} may be written in 
the form  
\begin{equation} \label{ZEqnNum577970} D_{x^{i} }^{\alpha } L-d_{t}^{\alpha } (D_{y^{i(\alpha )} }^{\alpha } L)-y^{j(2\alpha )} D_{y^{j(\alpha )} }^{\alpha } (D_{y^{i(\alpha )} }^{\alpha } L)=0, \end{equation} 
for $i=\overline{1, n}$. The equations 
\eqref{ZEqnNum578225} may be written as 
\begin{equation} \label{56)} \frac{\partial L}{\partial x^{i} } -d_{t}^{\alpha } \left(\frac{\partial L}{\partial y^{i(\alpha )} } \right)-y^{j(2\alpha )} \frac{\partial ^{2} L}{\partial y^{i(\alpha )} \partial y^{j(\alpha )} } =0, \end{equation} 
where $i=\overline{1, n}$. 
Let us denote by 
\begin{equation} \label{57)} g_{ij}^{\alpha } =D_{y^{i(\alpha )} }^{\alpha } (D_{y^{j(\alpha )} }^{\alpha } L), \end{equation} 
and by 
$\left(\mathop{g^{ik} }\limits^{\alpha } \right)=\left(g_{ij}^{\alpha } \right)^{-1} $, 
if $\det (g_{ij}^{\alpha } )\ne 0$. 
From \eqref{ZEqnNum577970} and from 
Proposition 5, we get the fractional 
field \textit{FODE} 
$\mathop{\Gamma _{L} }\limits^{\alpha }$ 
associated to $L$. \\ 
Let $c:t\in [0,1]\to (x^{i} (t))\subset U$ be a 
parameterized curve. 
The extension of $c$ to 
$J^{\alpha k} ({\mathbb R},\; M)$ 
is the curve 
$c^{\alpha k} :t\in [0,1]\to (t,x^{i} (t),y^{\alpha a} (t))\in J^{\alpha k} ({\mathbb R},\; M)$, 
$a=\overline{1,k}$. 
Let 
$L:J^{\alpha k} ({\mathbb R},\; M)\to {\mathbb R}$ 
be a Lagrange function. The action of $L$ along 
the curve $c^{\alpha k} $ is 
\begin{equation} \label{ZEqnNum533362} {\mathcal A}(c^{\alpha k} )=\int _{0}^{1}L(t,x(t),y^{\alpha a} (t))dt . 
\end{equation} 
Let 
$c_{\varepsilon } :t\in [0,1]\to (x^{i} (t,\varepsilon ))\in M$ 
be a family of curves, where the absolute value of 
$\varepsilon $ is sufficiently small so that 
$Imc_{\varepsilon } \subset U\subset M$, 
$c_{0} (t)=c(t)$, 
$D_{\varepsilon }^{\alpha } c(\varepsilon )\left|_{\varepsilon =0} \right. =D_{\varepsilon }^{\alpha } c(\varepsilon )\left|_{\varepsilon =1} \right. =0$. 
The action of $L$ on the curve 
$c_{\varepsilon } $ is given by 
\begin{equation} \label{ZEqnNum472103} {\mathcal A}(c_{\varepsilon }^{\alpha k} )=\int _{0}^{1}L(t,x(t,\varepsilon ),y^{\alpha a} (t,\varepsilon ))dt  \end{equation} 
where 
$y^{i(\alpha a)} (t,\varepsilon )=\frac{1}{\Gamma (1+\alpha a)} D_{t}^{\alpha a} x^{i} (t,\varepsilon )$, 
$a=\overline{1,k}$. 
The action \eqref{ZEqnNum472103} has a fractional 
extremal value if 
\begin{equation} \label{60)} D_{\varepsilon }^{\alpha } ({\mathcal A}(c_{\varepsilon }^{\alpha k} ))\left|_{\varepsilon =0} \right. =0. \end{equation} 
The action \eqref{ZEqnNum472103} has an extremal value if 
\begin{equation} \label{61)} D_{\varepsilon }^{1} ({\mathcal A}(c_{\varepsilon }^{\alpha k} ))\left|_{\varepsilon =0} \right. =0. \end{equation} 
\begin{proposition} 
a) A necessary condition for the action 
\eqref{ZEqnNum533362} to reach a fractional 
extremal value is that $c(t)$ satisfies 
the fractional Euler-Lagrange equations 
\begin{equation} \label{eq62} D_{x^{i} }^{\alpha } L+\sum _{a=1}^{k}(-1)^{a} d_{t}^{\alpha a} (D_{y^{i(\alpha a)} }^{\alpha } L)=0 , \end{equation} 
where 
\begin{equation} \label{63)} d_{t}^{\alpha a} =D_{t}^{\alpha } +y^{i(\alpha )} D_{x^{i} }^{\alpha } +y^{i(2\alpha )} D_{y^{i(\alpha )} }^{\alpha } +...+y^{i(\alpha a)} D_{y^{i(\alpha (a-1))} }^{\alpha } , \end{equation} 
and $i=\overline{1, n}$. \\ 
b) A necessary condition that the action 
\eqref{ZEqnNum533362} reaches an extremal 
value is that $c(t)$ satisfies the 
Euler-Lagrange equations 
\begin{equation} \label{eq64} \frac{\partial L}{\partial x^{i} } +\sum _{a=1}^{k}(-1)^{a} d_{t}^{a} (D_{y^{i(\alpha a)} }^{\alpha } L)=0 , \end{equation} 
where 
\begin{equation} \label{ZEqnNum420605} d_{t}^{a} =D_{t}^{1} +y^{i(\alpha )} D_{x^{i} }^{1} +...+y^{i(\alpha a)} D_{y^{i(\alpha (a-1))} }^{1} . \end{equation} 
\end{proposition} 
\textbf{Example.} 
Consider the fractional differential equation 
\begin{equation} \label{eq66} \begin{array}{l} {\frac{c\Gamma (1+\gamma )}{\Gamma (1+\gamma -\alpha )} x^{\gamma -\alpha } (t)f(t)+a_{1} \Gamma (1+2\alpha )y^{(2\alpha )} +} \\ {a_{2} \Gamma (1+3\alpha )y^{(3\alpha )} =0.} \end{array} \end{equation} 
The equation \eqref{eq66} is the fractional 
Euler-Lagrange equation \eqref{eq62} 
for the function 
\begin{equation*} \begin{array}{l} {L=\frac{c}{1+\gamma -\alpha } x^{\gamma } -a_{1} \Gamma (1+2\alpha )(y^{\alpha } )^{\alpha } +} \\ {a_{2} \Gamma (1+3\alpha )(y^{2\alpha } )^{\alpha } .} \end{array} \end{equation*} 
The equation \eqref{eq66} is the fractional 
Euler-Lagrange equation \eqref{eq64} for the 
function 
\begin{equation*}  L=\frac{c\Gamma (1+\gamma )x^{\gamma -\alpha +1} }{\Gamma (1+\gamma -\alpha )^{(1+\gamma -\alpha )} } f-\frac{a_{1} }{2} \Gamma (1+2\alpha )(y^{\alpha } )^{2} +\frac{a_{2} }{2} \Gamma (1+3\alpha )(y^{2\alpha } )^{2} . \end{equation*} 
\section{Examples and applications} 
1. \textbf{The nonhomogeneous Bagley-Torvik equation} \\ 
\indent The dynamics of a flat rigid body embedded 
in a Newton fluid is described by the equation 
\begin{equation} \label{ZEqnNum785447} aD_{t}^{2}x(t)+bD_{t}^{3/2} x(t)+cx(t)-f(t)=0, \end{equation} 
where $a,b,c\in {\mathbb R}$ and the initial conditions 
are $x(0)=0$,  $D_{t}^{1}x(0)=0$. 
The equation \eqref{ZEqnNum785447} is a fractional 
differential equation on the bundle 
$J^{\alpha } ({\mathbb R},\; {\mathbb R})$ 
for $\alpha =\frac{1}{4} $. Indeed, let's consider 
the fractional differential equation 
\begin{equation} \label{ZEqnNum769966} aD_{t}^{8\alpha } x(t)+bD_{t}^{6\alpha } x(t)+cx(t)-f(t)=0, \end{equation} 
with $\alpha >0$. For $\alpha =\frac{1}{4}$ 
the equation \eqref{ZEqnNum769966} reduces to 
\eqref{ZEqnNum785447}. With the notations 
\eqref{ZEqnNum283756}, the equation 
\eqref{ZEqnNum769966} becomes 
\begin{equation} \label{68)} a\Gamma (1+8\alpha )y^{(8\alpha )} (t)+b\Gamma (1+6\alpha )y^{(6\alpha )} (t)+cx(t)-f(t)=0. \end{equation} 
On the bundle $J^{4\alpha } ({\mathbb R},\; {\mathbb R})$ 
let us consider the Lagrange function 
\begin{equation} \label{ZEqnNum983897} \begin{array}{l} {L(t,x,y^{(3\alpha )} ,y^{(4\alpha )} )=\frac{1}{2} cx^{2} -fx-\frac{b}{2} \Gamma (1+6\alpha )(y^{(3\alpha )} )^{2} +} \\ {\frac{a}{2} \Gamma (1+8\alpha )(y^{(4\alpha )} )^{2} .} \end{array} \end{equation} 
Using the relation \eqref{ZEqnNum420605}, 
the Euler-Lagrange equation for \eqref{ZEqnNum983897} 
is 
\begin{equation} \label{70)} \begin{array}{l} {D_{x}^{1} L-D_{t}^{3\alpha } (D_{y^{(3\alpha )} }^{1} L)+D_{t}^{4\alpha } (D_{y^{(4\alpha )} }^{1} L)=} \\ {cx-f+b\Gamma (1+6\alpha )D_{t}^{3\alpha } y^{(3\alpha )} +a\Gamma (1+8\alpha )D_{t}^{4\alpha } y^{(4\alpha )} =} \\ {cx-f+b\Gamma (1+6\alpha )y^{(6\alpha )} +a\Gamma (1+8\alpha )y^{(8\alpha )} =0.} \end{array} \end{equation} 
\begin{proposition} 
The equation \eqref{ZEqnNum785447} represents 
the Euler-Lagrange equation on the bundle 
$J^{4\alpha } ({\mathbb R},\; {\mathbb R})$ for 
$\alpha =\frac{1}{4} $, with the Lagrange function 
given by 
\begin{equation} \label{71)} \begin{array}{l} {L(t,x,y^{(3/2)} ,y^{(2)} )=\frac{1}{2} cx^{2} -fx-\frac{b}{2} \Gamma ({5\mathord{\left/ {\vphantom {5 2)}} \right. \kern-\nulldelimiterspace} 2)} (y^{(3/2)} )^{2} +} \\ {\frac{a}{2} \Gamma (3)(y^{(2)} )^{2} .} \end{array} \end{equation} 
\end{proposition} 
2. \textbf{Differential equations of order 
one, two and three which admit fractional 
Lagrangians} \\ 
The following differential equations 
don't have classical Lagrangians such that 
the Euler-Lagrange equation represents 
the given equation: 
\begin{equation} \label{ZEqnNum860560} \dot{x}(t)+V_{1} (t,x)=0,\; \; V_{1} (t,x)=\frac{\partial U_{1} (t,x)}{\partial x} , \end{equation} 
\begin{equation} \label{ZEqnNum277412} \ddot{x}(t)+a_{1} \dot{x}(t)+V_{2} (t,x)=0,\; \; V_{2} (t,x)=\frac{\partial U_{2} (t,x)}{\partial x} , \end{equation} 
\begin{equation} \label{ZEqnNum426478} \dddot{x}(t)+a_{2} \ddot{x}(t)+a_{1} \dot{x}(t)+V_{3} (t,x)=0,\; \; V_{3} (t,x)=\frac{\partial U_{3} (t,x)}{\partial x} . 
\end{equation} 
Let us associate the fractional equations from below 
to the equations \eqref{ZEqnNum860560}, 
\eqref{ZEqnNum277412} and 
\eqref{ZEqnNum426478}, respectively: 
\begin{equation} \label{75)} D_{t}^{2\alpha } x(t)+V_{1} (t,x)=0, \end{equation} 
\begin{equation} \label{76)} D_{t}^{4\alpha } x(t)+a_{1} D_{t}^{2\alpha } x(t)+V_{2} (t,x)=0, \end{equation} 
\begin{equation} \label{77)} D_{t}^{6\alpha } x(t)+a_{2} D_{t}^{4\alpha } x(t)+a_{1} D_{t}^{2\alpha } x(t)+V_{3} (t,x)=0. \end{equation} 
\begin{proposition} 
a) Let 
$J^{\alpha}({\mathbb R},\;  {\mathbb R})\to {\mathbb R}$ 
be the fractional bundle and consider 
$L:J^{\alpha}({\mathbb R},\; {\mathbb R})\to {\mathbb R}$ 
given by 
\begin{equation} \label{ZEqnNum972723} L(t,x,y^{(\alpha )} )=U_{1} (t,x)-\frac{1}{2} \Gamma (1+2\alpha )(y^{\alpha } )^{2} . \end{equation} 
The Euler-Lagrange equation of \eqref{ZEqnNum972723} is 
\begin{equation} \label{79)} \begin{array}{l} {\frac{\partial L}{\partial x} -D_{t}^{\alpha } \left(\frac{\partial L}{\partial y^{\alpha } } \right)=\frac{\partial U_{1} (t,x)}{\partial x} +\Gamma (1+2\alpha )y^{(2\alpha )} =} \\ {V_{1} (t,x)+D_{t}^{2\alpha } x(t)=0.} \end{array} \end{equation} 
b) Let 
$J^{2\alpha } ({\mathbb R},\; {\mathbb R})\to {\mathbb R}$ 
be the fractional bundle and the Lagrangian 
$L:J^{2\alpha } ({\mathbb R},\; {\mathbb R})\to {\mathbb R}$ 
given by 
\begin{equation} \label{ZEqnNum400853} \begin{array}{l} {L(t,x,y^{(\alpha )} ,y^{(2\alpha )} )=U_{2} (t,x)-\frac{1}{2} a_{1} \Gamma (1+2\alpha )(y^{\alpha } )^{2} +} \\ {\frac{1}{2} \Gamma (1+4\alpha )(y^{(2\alpha )} )^{2} .} \end{array} \end{equation} 
The Euler-Lagrange equation of \eqref{ZEqnNum400853} 
is 
\begin{equation} \label{81)} \begin{array}{l} {\frac{\partial L}{\partial x} -D_{t}^{\alpha } \left(\frac{\partial L}{\partial y^{\alpha } } \right)+D_{t}^{2\alpha } \left(\frac{\partial L}{\partial y^{2\alpha } } \right)=} \\ {V_{2} (t,x)+a_{1} \Gamma (1+2\alpha )y^{(2\alpha )} +} \\ {a_{2} \Gamma (1+4\alpha )y^{(4\alpha )} =} \\ {V_{2} (t,x)+a_{1} D_{t}^{2\alpha } x(t)+D_{t}^{4\alpha } x(t)=0.} \end{array} \end{equation} 
c) Let 
$J^{3\alpha } ({\mathbb R},\; {\mathbb R})\to {\mathbb R}$ 
be the fractional bundle and 
$L:J^{3\alpha } ({\mathbb R},\; {\mathbb R})\to {\mathbb R}$ 
given by 
\begin{equation} \label{ZEqnNum619704} \begin{array}{l} {L(t,x,y^{(\alpha )} ,y^{(2\alpha )} ,y^{(3\alpha )} )=V_{3} (t,x)-\frac{a_{1} }{2} \Gamma (1+2\alpha )(y^{(\alpha )} )^{2} +} \\ {\frac{a_{2} }{2} \Gamma (1+4\alpha )(y^{(2\alpha )} )^{2} -\frac{1}{2} \Gamma (1+6\alpha )(y^{(3\alpha )} )^{2} .} \end{array} \end{equation} 
The Euler-Lagrange equation of \eqref{ZEqnNum619704} 
is 
\begin{equation} \label{83)} \begin{array}{l} {\frac{\partial L}{\partial x} -D_{t}^{\alpha } \left(\frac{\partial L}{\partial y^{\alpha } } \right)+D_{t}^{2\alpha } \left(\frac{\partial L}{\partial y^{(2\alpha )} } \right)-D_{t}^{3\alpha } \left(\frac{\partial L}{\partial y^{(3\alpha )} } \right)=V_{3} (t,x)+} \\ {a_{1} \Gamma (1+2\alpha )y^{(2\alpha )} +a_{2} \Gamma (1+4\alpha )y^{(4\alpha )} +\Gamma (1+6\alpha )y^{(6\alpha )} =} \\ {V_{3} (t,x)+a_{1} D_{t}^{2\alpha } x(t)+a_{2} D_{t}^{4\alpha } x(t)+D_{t}^{6\alpha } x(t)=0.} \end{array} \end{equation} 
d) For $\alpha =\frac{1}{2} $ we obtain the fractional 
Lagrangians that describe the equations 
\eqref{ZEqnNum860560}, \eqref{ZEqnNum277412}, 
\eqref{ZEqnNum426478}, respectively 
\begin{equation} \label{84)} \begin{array}{l} {L(t,x,y^{(1/2)} )=U_{1} (t,x)-\frac{1}{2} \Gamma (2)(y^{(1/2)} )^{2} } \\ {L(t,x,y^{(1/2)} ,y^{(1)} )=U_{2} (t,x)-\frac{1}{2} a_{1} \Gamma (2)(y^{(1/2)} )^{2} +\frac{1}{2} \Gamma (3)(y^{(1)} )^{2} } \\ {L(t,x,y^{(1/2)} ,y^{(1)} ,y^{(3/2)} )=U_{3} (t,x)-\frac{a_{1} }{2} \Gamma (2)(y^{(1/2)} )^{2} +} \\ {\frac{a_{2} }{2} \Gamma (2)(y^{(1)} )^{2} -\frac{1}{2} \Gamma (4)(y^{(3/2)} )^{2} .} \end{array} \end{equation} 
\end{proposition} 
In the category of the equations 
\eqref{ZEqnNum277412} and \eqref{ZEqnNum426478} 
there are: \\ 
a) the nonhomogeneous classical friction equation 
\begin{equation} \label{85)} m\ddot{x}(t)+\gamma \dot{x}(t)-\frac{\partial U(t,x)}{\partial x} =0, \end{equation} 
b) the nonhomogeneous model of Phillips \cite{Lorenz} 
\begin{equation} \label{86)} \ddot{x}(t)+a_{1} \dot{x}(t)+b_{1} x(t)+f(t)=0, \end{equation} 
c) the nonhomogeneous business cycle with 
innovation \cite{Lorenz} 
\begin{equation} \label{87)} \dddot{y}(t)+a_{2} \ddot{y}(t)+a_{1} \dot{y}(t)+b_{1} x(t)+f(t)=0. 
\end{equation} 
\section*{Conclusions} 
The paper presents the main differentiable structures 
on $J^{\alpha } ({\mathbb R},\; M)$, in order to 
describe fractional differential equations and 
ordinary differential equations, using Lagrange functions 
defined on $J^{\alpha } ({\mathbb R},\; M)$. \\ 
\indent With the help of the methods shown, 
there may be analyzed other models, such as 
those found in \cite{Caputo} and \cite{broasca}.    
 
\end{document}